\documentstyle[12pt]{article}
\begin{document}
\large
\newtheorem{prop}{Proposition}
\newtheorem{cor}{Corollary}
\newtheorem{teo}{Theorem}
\title{Analysis and Probability over Infinite Extensions of a Local Field}
\author{Anatoly N. Kochubei\\ \footnotesize Institute of Mathematics,\\
\footnotesize Ukrainian National Academy of Sciences,\\
\footnotesize Tereshchenkivska 3, Kiev, 252601 Ukraine}
\date{}
\maketitle
\begin{abstract}
We consider an infinite extension $K$ of a local field of zero
characteristic which is a union of an increasing sequence of finite
extensions. $K$ is equipped with an inductive limit topology;
its conjugate $\overline{K}$ is a completion of $K$ with respect to
a topology given by certain explicitly written
seminorms. We construct and study a Gaussian measure, a Fourier
transform, a fractional differentiation operator and a cadlag
Markov process on $\overline{K}$. If we deal with Galois extensions
then all these objects are Galois-invariant.


\end{abstract}
\section{Introduction}

Much of the recent activity in $p$-adic analysis has been concentrated
on pseudo-differential operators acting on spaces of complex-valued
functions and distributions over local fields [VVZ, Ha2, K1-K4],
potential theory and stochastic processes on local fields [AK, Bl,
E2, E4, Ha1, K1, K5, KV, Va, Y].

All those results, as well as their counterparts in real analysis,
use in a very strong way the local compactness of the underlying field;
in particular they rely on existence of the Haar measure. On the
other hand, the calculus over
infinite-dimensional vector spaces where no invariant measure can
exist is among well-established topics of real analysis; see [AH,
BK, DF, Hi, Kuo, VTC].
Some general principles of such a calculus for the $p$-adic situation
were introduced by Evans [E1]; the first constructive results were
obtained by Madrecki [Ma] and Satoh [Sa].

The well-known Minlos theorem (see e.g. [DF]) states that
every continuous positive-definite function on a real nuclear
locally convex topological vector space $X$ is the Fourier transform
of some Radon measure  on the  conjugate
space $X^\ast $ equipped
with the $\ast $-weak topology. It was shown by Madrecki [Ma] that
a similar result in the $p$-adic case holds without any nuclearity
assumptions. This agrees with the fact [PT] that any locally
convex space over a $p$-adic local field is nuclear in the sense
of Grothendieck.

Satoh's paper [Sa] is devoted to a $p$-adic version of the theory of
abstract Wiener spaces (see [Kuo]). This leads to a construction of
Wiener-type measures on some non-Archimedean Banach spaces; in
particular, a certain space of power series has been considered.

Note that while a large part of the real infinite-dimensional analysis is
devoted to the study of measures, function spaces and operators over a
Hilbert space, there is no clear counterpart of a Hilbert space in the
$p$-adic case. However there are other infinite-dimensional spaces over
$p$-adics which are of purely arithmetical nature and constitute a natural
arena for developing analysis. These are infinite extensions of local fields.

In this paper we consider an infinite extension $K$ of a local field $k$,
$\mbox{char }k=0$, which is a union of an increasing sequence
$k=K_1\subset \ldots \subset K_n\subset \ldots $ of finite extensions.
The field $K$ is a topological vector space over $k$ with the inductive
limit topology. Its conjugate $\overline{K}$ is a completion of $K$ with
respect to a certain topology defined in arithmetical terms. We construct a
Radon measure on $\overline{K}$ which is Gaussian in the sense of Evans [E1]
and possesses some (partial) invariance properties. A version of the
Fourier-Wiener transform is introduced over $\overline{K}$, and Fourier
images of certain test functions are described. This allows to define and
study a pseudo-differential operator over $\overline{K}$ similar to the
fractional differentiation operator over a local field extensively
investigated in [VVZ, Ha1, K1, K4]. This operator is proved to be a
generator of
a Markov process on $\overline{K}$, an analogue of the symmetric stable
process. If we deal with Galois extensions then all these objects are
invariant with respect to the Galois group of the extension $K/k$.

The class of extensions under consideration includes the maximal
unramified extension (here $K_n/k$ is the unramified extension of
the degree $n!,\ n=2,3,\ldots $). Some other examples can be found in [FV].

\section{An Infinite Extension as a Topological Vector Space}

1. Let $k$ be a non-Archimedean local field of zero characteristic
(thus $k$ is a finite extension of the field of $p$-adic numbers;
for the basics about local fields see [CF, FV, W]). Consider an
increasing sequence of its finite extensions
\[
k=K_1\subset K_2\subset \ldots \subset K_n\subset \ldots .
\]
For each $n=1,2,\ldots $ we denote by $|\cdot |_n$ a normalized absolute
value in $K_n$; if $N_n:\ K_n\to k$ is a norm mapping then $|x|_n=
|N_n(x)|_1$. We shall use also a (non-normalized) absolute value $\|x\|=
|x|_n^{1/m_n}\ (x\in K_n)$ which coincides with $|x|_1$ for $x\in k$ being
defined correctly on the infinite extension $K=\bigcup \limits _{n=1}^\infty K_n$.
Here $m_n$ is the degree of the extension $K_n/k$.

Let us consider, for each $n$, a mapping $T_n:\ K\to K_n$ defined as follows.
If $x\in K_\nu ,\ \nu >n$, put
\[
T_n(x)=\frac{m_n}{m_\nu }\mbox{Tr}_{K_\nu /K_n}(x)
\]
where $\mbox{Tr}_{K_\nu /K_n}:\ K_\nu \to K_n$ is the trace mapping
(note that the degree of the extension $K_\nu /K_n$ equals $m_\nu /m_n$ ).
If $x\in K_n$ then, by definition, $T_n(x)=x$. It follows from the
transitivity property of traces that $T_n$ is well-defined and $T_n\circ
T_\nu =T_n$ for $\nu >n$. Below we shall often write $T$ instead of $T_1$.

Note that the equality $T_n(x)=0$ for all $n$ implies $x=0$. Indeed, if
$x\in K$ then $x\in K_l$ for some $l$, and $T_l(x)=x$.

The field $K$ may be considered as a separable topological vector space
over $k$ with the inductive limit topology (see [T] for the theory of
topological vector spaces over local fields).

Denote by $\overline{K}$ the projective limit of the sequence $\{K_n\}$
with respect to the mappings $\{T_n\}$. It is defined as a subset of the
direct product $\prod \limits _{n=1}^\infty K_n$ consisting of those
$x=(x_1,\ldots ,x_n,\ldots ),\ x_n\in K_n$, for which $x_n=T_n(x_\nu )$
if $\nu >n$. The topology in $\overline{K}$ is introduced by seminorms
\[
\|x\|_n=\|x_n\|,\ \ n=1,2,\ldots .
\]
Both spaces $K$ and $\overline{K}$ are complete and reflexive [T, Mo].

An element $x\in K$ may be identified with $(x_1,\ldots ,x_n,\ldots )\in
\overline{K}$ where $x_n=T_n(x)$. Thus $K$ may be viewed as a subset of
$\overline{K}$. the mappings $T_n$ can be extended to linear continuous
mappings from $\overline{K}$ to $K_n$ by setting $T_n(x)=x_n$ for any
$x=(x_1,\ldots ,x_n,\ldots )\in \overline{K}$.

\begin{prop}
The topological vector spaces $K$ and $\overline{K}$ are the strong duals
of each other, with the pairing defined by the formula
\[
\langle x,y\rangle =T(xy_n)
\]
where $x\in K_n\subset K,\ y=(y_1,\ldots ,y_n,\ldots )\in \overline{K},\
y_n\in K_n$. As a subset of $\overline{K},\ K$ is dense in $\overline{K}$.
\end{prop}

{\it Proof.} Each field $K_n$ may be identified as a topological vector
space over $k$ with its dual $K_n^\ast $ -- if $\varphi \in K_n^\ast $
then there exists a unique element $a_\varphi \in K_n$ such that
\[
\varphi (x)=T(a_\varphi x),\ \ x\in K_n,
\]
and $\varphi $ may be identified with $a_\varphi $.

Let $\nu >n$. The imbedding $\tau _{n\nu }:\ K_n\to K_\nu $
determines (due to the above identification) the adjoint
mapping $\tau _{n\nu }^\ast :\ K_\nu \to K_n$. We shall show that
$\tau _{n\nu }^{\ast }=T_n|K_\nu $.

Indeed, if $x\in K_n,\ \psi \in K_\nu ^\ast $ then
\[
\psi (\tau _{n\nu }(x))=\varphi (x)
\]
where $\varphi =\tau _{n\nu }^\ast (\psi )\in K_n^\ast$.
This means that for any $x\in K_n$
\[
\frac{1}{m_\nu }\mbox{Tr}_{K_{\nu }/k}(a_\psi x)=\frac{1}{m_n}
\mbox{Tr}_{K_n/k}(a_\varphi x),
\]
$a_\varphi \in K_n,\ a_\psi \in K_\nu $. On the other hand,
\[
\mbox{Tr}_{K_{\nu }/k}(a_\psi x)=\mbox{Tr}_{K_n/k}(\mbox{Tr}_
{K_{\nu }/K_n}(a_\psi x))=\mbox{Tr}_{K_n/k}(x\mbox{Tr}_
{K_{\nu }/K_n}(a_\psi ))
\]
whence $a_\varphi =\frac{m_n}{m_\nu }\mbox{Tr}_{K_{\nu }/K_n}(a_\psi )
=T_n(a_\psi )$.

Now the duality result follows from the Komatsu-Morita duality theorem [Mo].

Let $y=(y_1,\ldots ,y_n,\ldots )\in \overline{K}$. Consider a
sequence $\{y^{(n)}\}\subset \overline{K}$ of elements $y^{(n)}
=y_n\in K$ viewed as lying in $\overline{K}$. For each
$j=1,2,\ldots $ we have
\[
T_j(y_n)=y_j,\ \ n\ge j.
\]
Hence $\|y^{(n)}-y\|_j\to 0,\ \ n\to \infty $, for each $j$,
so that $K$ is dense in $\overline{K}.\ \ \ \ \Box $

\begin{cor}
The space $\overline{K}$ is separable.
\end{cor}
{\it Proof.}  One can take a countable dense subset in
each $K_n$; their union is countable and dense in $K$ with
respect to the topology induced from $\overline{K}$, thus it is
dense in $\overline{K}.\ \ \ \ \Box $

2. Suppose that $K_n/k$ is a Galois extension for each $n$
with the Galois group $G_n$. Then $K/k$ is a Galois extension [Bo],
its Galois group $G$ is a projective limit of the groups $G_n$ [K].
The group $G$ consists of sequences
$g=(g_1,\ldots ,g_n,\ldots ),\ g_n\in G_n$, such that the
restriction of an automorphism $g_\nu ,\ \nu >n$, to $K_n$
coincides with $g_n$.
The group operations in $G$ are defined component-wise;
if $x\in K$ then $x\in K_n$ for some $n$, and $g(x)=g_n(x)$.

Let us consider elements $x,y\in K$. Suppose that $y\in K_n$
for some $n$. Then for any $g\in G$
\[
T(g(x)y)=T(g_n(x)y)=m_n^{-1}\mbox{Tr}_{K_n/k}(g_n(x)y)
=m_n^{-1}\mbox{Tr}_{K_n/k}(xg_n^{-1}(y))
\]
due to the $G$-invariance of the trace. Hence,
\begin{equation}
T(g(x)y)=T(xg^{-1}(y)),\ \ x,y\in K.
\end{equation}

The relation (1) shows that the action of $G$ on $K$ is continuous
with respect to the topology induced on $K$ by the
$\ast $-weak topology of $\overline{K}$. Since $K$ is strongly
(and, of course, $\ast $-weakly) dense in $\overline{K}$, we
come to the following result.

\begin{prop}
The action of the Galois group $G$ on $K$ can be extended by
continuity to its $\ast $-weakly continuous action on $\overline{K}$.
\end{prop}

\section{Measure and Integration}

1. Let us start from considering various $\sigma $-algebras of
subsets in $\overline{K}$.

Denote by ${\cal B}(\overline{K})$ the Borel $\sigma $-algebra
over $\overline{K}$ with respect to the strong topology of the
(non-Archimedean) Frechet space $\overline{K}$; the Borel
$\sigma $-algebra with respect to the $\ast $-weak topology
will be denoted by ${\cal B}(\overline{K},K)$. We shall consider
also the smallest  $\sigma $-algebra ${\cal A}(\overline{K})$
for which all the linear continuous functionals $\overline{K}\to k$
are measurable.

\begin{prop}
The above  $\sigma $-algebras coincide:
\[
{\cal A}(\overline{K})={\cal B}(\overline{K},K)={\cal B}(\overline{K}).
\]
\end{prop}
{\it Proof.}  It is clear that
\[
{\cal A}(\overline{K})\subset {\cal B}(\overline{K},K)\subset
{\cal B}(\overline{K}).
\]
Hence it is sufficient to prove that ${\cal A}(\overline{K})
={\cal B}(\overline{K})$ or (by definition of the topology in
$\overline{K}$) that the set
\[
V_n=\{x\in \overline{K}:\ \ \|x\|_n\le 1\}
\]
belongs to ${\cal A}(\overline{K})$.

Consider its polar
\[
V_n^0=\{\xi \in K:\ \ |\langle \xi ,x\rangle |_1\le 1\ \ \forall
x\in V_n\}.
\]
By the bipolar theorem ([T], Theorem 4.14)
\[
V_n=\{x\in \overline{K}:\ \ \sup \limits _{\xi \in V_n^0}\
|\langle \xi ,x\rangle |_1\le 1\}.
\]

On the other hand, $V_n^0$ is equicontinuous [T] (here we look at
$K$ as the dual of $\overline{K}$). Since $\overline{K}$ is separable
we find (just as in the Archimedean case, see [Sch]) that the $\ast $-weak
topology on $V_n^0$ is metrisable. Now being $\ast $-weakly
compact [T] and metrisable, the set $V_n^0$ is separable in
the $\ast $-weak topology. This implies the existence of such a
countable subset $M\subset V_n^0$ that
\begin{displaymath}
V_n=\{x\in \overline{K}:\ \ \sup \limits _{\xi \in M}\ |\langle
\xi ,x\rangle |_1\le 1\}
=\bigcap_{\xi \in M}\{x\in \overline{K}:\ \ |\langle \xi ,x\rangle |_1\le 1\}.
\end{displaymath}
Thus $V_n\in {\cal A}(\overline{K})$.\ \ \ $\Box $.

Similar results in the framework of real analysis can be found in [VTC].

2. Let us recall (in a form convenient for subsequent use) some
results from harmonic analysis on local fields [GGP, ST].

Denote by $q_n$ the residue field cardinality for the field $K_n$,
and by $d_n$ the exponent of the different for the extension $K_n/k$.
The Haar measure on the additive group of the field $K_n$ will be
denoted $dx$ (the dependence on $n$ is not shown; however this will
not lead to any confusion). This measure will be normalized by the relation
\[
\int _{\|x\|\le 1}dx=1.
\]
Using the relation $d(ax)=|a|_n\,dx,\ a\in K_n$, we find that
\[
\int _{|x|_n\le q_n^\nu }dx=q_n^\nu ;\ \ \int _{|x|=q_n^\nu }dx
=(1-q_n^{-1})q_n^\nu ,\ \ \nu \in \mbox{{\bf Z}}.
\]

Let $\chi $ be a rank zero additive character on $k$. Then for
any $a\in K_n,\ \nu \in \mbox{{\bf Z}}$
\begin{equation}
\int _{|x|_n\le q_n^\nu }\chi \circ \mbox{Tr}_{K_n/k}(ax)\,dx
=\left\{ \begin{array}{rl}
q_n^\nu ,& \mbox{if } |a|_n\le q^{d_n-\nu }_n\\
0, & \mbox{if } |a|_n> q^{d_n-\nu }_n
\end{array}\right.
\end{equation}
\begin{equation}
\int _{|x|_n=q_n^\nu }\chi \circ \mbox{Tr}_{K_n/k}(ax)\,dx
=\left\{ \begin{array}{rl}
(1-q_n^{-1})q_n^\nu ,& \mbox{if } |a|_n\le q^{d_n-\nu }_n\\
-q_n^{\nu -1},& \mbox{if } |a|_n=q^{d_n-\nu +1}_n\\
0, & \mbox{if } |a|_n\ge q^{d_n-\nu +2}_n
\end{array}\right.
\end{equation}

The Fourier transform (defined initially on $L_1(K_n)$, with an
isometric extension onto $L_2(K_n)$) will be defined as
\begin{equation}
\widetilde{f}(\xi )=q_n^{-d_n/2}\int _{K_n}\chi \circ \mbox{Tr}_
{K_n/k}(x\xi)f(x)\,dx\,,
\end{equation}
with the inversion formula
\begin{equation}
f(x)=q_n^{-d_n/2}\int _{K_n}\chi \circ \mbox{Tr}_{K_n/k}(-x\xi)
\widetilde{f}(\xi )\,d\xi\,,
\end{equation}
valid,for example, for a locally constant complex valued
function $f$ with a compact support. The class ${\cal D}(K_n)$
of such functions is invariant under the Fourier transform.

In particular, if
\[
\Omega _n(x)=\left\{ \begin{array}{rl}
1 ,& \mbox{if } \|x\|\le 1\\
0, & \mbox{if } \|x\|>1
\end{array}\right.
\]
$(x\in K_n)$ then by (2) and (4)
\begin{equation}
\widetilde{\Omega _n}(\xi )=\left\{ \begin{array}{rl}
q_n^{-d_n/2} ,& \mbox{if } \|\xi \|\le q_n^{d_n/m_n}\\
0, & \mbox{if } \|\xi \|>q_n^{d_n/m_n}
\end{array}\right.
\end{equation}
$(\xi \in K_n)$.

Note that a cutoff of the Haar measure concentrated on some
ball $\{x\in K_n:\ |x|_n\le q_n^\nu \}$ was found by Evans [E1]
to be the only natural $p$-adic counterpart of the Gaussian measure.

3.Let us consider a function $\Omega (x),\ x\in K$, given by the formula
\[
\Omega (x)=\Omega _n(x),\ \ x\in K_n.
\]
This function (which is clearly well-defined) is continuous on $K$
by the definition of the inductive limit topology. It follows from (6)
and (5) that $\Omega _n$ is a positive definite on each $K_n$.
Hence, $\Omega $ is a positive definite function on $K$.

Using Madrecki's theorem [Ma] and Propositions 1,3 we establish the
existence of a unique Radon probability measure $\mu $ on
${\cal B}(\overline{K})$ such that
\begin{equation}
\Omega (a)=\int _{\overline{K}}\chi (\langle a,x\rangle )\,d\mu
(x),\ \ a\in K.
\end{equation}

We begin our study of properties of the measure $\mu $ with
integrating a "cylindrical" function of the form $f(x)=\varphi
(T_n(x))$ where $\varphi \in {\cal D}(K_n)$. The set of all such
functions will be denoted by ${\cal D}(\overline{K})$. Note that
${\cal D}(\overline{K})$ is a linear set - if $f(x)
=\varphi (T_n(x))+\psi (T_\nu (x)),\ \ \nu >n$, then $f(x)=(\varphi \circ
T_n+\psi )(T_\nu (x))$.

The function $\varphi $ can be written as
\[
\varphi (y)=\int _{K_n}\chi (T(ay))\psi (a)\,da,\ \ \psi \in {\cal D}(K_n),
\]
where $\psi $ is the inverse Fourier transform of the function
$q_n^{-d_n/2}\varphi (m_ny)$. By (7) and Fubini theorem, using the
relation $T\circ T_n=T$ we find that
\[
\int _{\overline{K}}f(x)\,d\mu (x)
=\int _{a\in K_n:\ \|a\|\le 1}\psi (a)\,da.
\]
In view of (6) the Plancherel formula implies
\[
\int _{\overline{K}}f(x)\,d\mu (x)=q_n^{-d_n}\int \limits_{y\in K_n:\
\|y\|\le q_n^{d_n/m_n}}\varphi (m_ny)\,dy
\]
whence
\begin{equation}
\int _{\overline{K}}f(x)\,d\mu (x)
=q_n^{-d_n}\|m_n\|^{-m_n}\int \limits_{z\in K_n:\
\|z\|\le q_n^{d_n/m_n}\|m_n\|}\varphi (z)\,dz.
\end{equation}

By an approximation argument, we can easily show now that the
formula (8) is valid for any function $f(x)=\varphi (T_n(x))$
with $\varphi \in L_1(K_n)$.

In the terminology of Evans [E1], Proposition 3 means that the
pair \linebreak $(\overline{K},{\cal B}(\overline{K}))$ is a measurable vector
space, so that there is a notion of a Gaussian measure on
$(\overline{K},{\cal B}(\overline{K}))$. The latter is understood
as a distribution of a Gaussian random variable which is defined as follows.

Let $X$ be a random variable with values in $k$ defined on some
probability space. $X$ is Gaussian if when $X_1,X_2$ are two
independent copies of $X$ and
$(\alpha _{11},\alpha _{12}),(\alpha _{21},\alpha _{22})\in k\oplus k$
are orthonormal then the pair $(X_1,X_2)$ has the same distribution
as the pair $(\alpha _{11}X_1+\alpha _{12}X_2,\alpha _{21}X_1+
\alpha _{22}X_2)$. Here orthogonality is understood in
non-Archimedean sense (see e.g. [NBB]).

\begin{prop}
$\mu $ is a Gaussian measure on $(\overline{K},{\cal B}(\overline{K}))$.
\end{prop}
{\it Proof.} Let $a\in K$, $B$ be a Borel subset of $k$,
$\varphi _B:\ k\to \mbox{{\bf R}}$ be an indicator of the set $B$.
Consider the function
\[
I(a)=\int _{\overline{K}}\varphi _B(\langle a,x\rangle )\,d\mu (x),
\ \ a\in K.
\]
By the results of Evans, it suffices to prove that
\begin{equation}
I(a)=q_1^m\int _{z\in B,|z|_1\le q_1^{-m}}dz
\end{equation}
for some $m\in \mbox{{\bf Z}}$.

Repeating the argument from the proof of (8) we find that
\[
I(a)=\int _{k}\varphi _B(z)\,dz\int _{|y|_1\le \|a\|^{-1}}\chi (-zy)\,dy.
\]
Let $a\in K_n,\ |a|_n=q_n^{N},\ N\in \mbox{{\bf Z}}$. Then $\|a\|=q_n^{N/m_n}
=q_1^{Nf_n/m_n}=q_1^{N/e_n}$ where $f_n,e_n$ are the modular degree
and the ramification index of the extension $K_n/k$.

By the identity (2)
\[
\int _{|y|_1\le \|a\|^{-1}}\chi (-zy)\,dy
= \int _{|y|_1\le q_1^{-N/e_n}}\chi (zy)\,dy
= \int _{|y|_1\le q_1^{[-N/e_n]}}\chi (zy)\,dy=
\]
\[
=\left\{ \begin{array}{rl}
q_1^{[-N/e_n]}, & \mbox{if } |z|_1\le q_1^{-[-N/e_n]}\\
0, & \mbox{if } |z|_1> q_1^{-[-N/e_n]},
\end{array}\right.
\]
and we have come to (9) with $m=[-N/e_n].\ \ \ \Box $

Several important properties of the measure $\mu $ are collected
in the following theorem.

\begin{teo}
(a) The measure $\mu $ is concentrated on the set
\[
S=\left\{ x\in \overline{K}:\ \ \|T_n(x)\|\le q_n^{d_n/m_n}\|m_n\|,\
n=1,2,\ldots \right\} .
\]
\newline (b) The measure $\mu $ is invariant with respect to
additive shifts by elements from $S$.\newline
(c) An additive shift by an element from $\overline{K}\setminus
S$ transforms $\mu $ into a measure orthogonal to $\mu $.
\newline (d) if $a\in K_n$ then
\begin{equation}
\mbox{ess sup }_{x\in \overline{K}}\ |
\langle a,x\rangle |_1=q_1^{-[-N/e_n]}
\end{equation}
where $|a|_n=q_n^N,\ N\in $ {\em {\bf Z}}. \newline
(e) The set ${\cal D}(\overline{K})$ is dense in $L_s(\overline{K})$
for any $s\in [1,\infty )$.
\end{teo}

{\bf Remarks.}  1. It may happen that $S\cap K=\{0\}$. For example,
if $K$ is the maximal unramified extension of $k$
then $d_n=0,\ \|m_n\|=|n!|_1\to 0$ for $n\to \infty $. On the other
hand, if $y\in K_r$ then $T_n(y)=y$ for $n\ge r$, and $y\in S$
implies $y=0$.

2. The equality (10) should be compared with the expression for
the $L_2$-norm of a linear functional on a real Hilbert space with a
Gaussian measure (see [DF]). As expected, in the non-Archimedean
case the $L_\infty $-norm has appeared. On the other hand, the
property (b) has no clear counterpart in the real
infinite-dimensional analysis. Meanwhile, such properties are
natural for the non-Archimedean case, in view of the general
results from [E3].

{\it Proof of Theorem 1.}  (a) Consider for each $n=1,2,\ldots $ the set
\[
S_n=\left\{ x\in \overline{K}:\ \ \|T_n(x)\|\le q_n^{d_n/m_n}\|m_n\|\right\} .
\]
Suppose that $x\in S_\nu ,\ \nu >n$. Then
\[
\|T_n(x)\|=\|m_n\|\cdot \|m_\nu \|^{-1}\|\mbox{Tr}_{K_\nu /K_n}(T_\nu (x))\|.
\]
Since $\|m_\nu ^{-1}T_\nu (x)\|\le q_\nu ^{d_\nu /m_\nu }$
we have
\[
|m_\nu ^{-1}T_\nu (x)|_n\le q_\nu ^{d_\nu }
\]
whence (see Chapter 8 in [W])
\[
|m_\nu ^{-1}\mbox{Tr}_{K_\nu /K_n}(T_\nu (x))|_n\le q_n^l
\]
where $l\in \mbox{{\bf Z}},\ e_{n\nu }(l-1)<d_\nu -d_{n\nu }\le e_{n\nu }l$,
$e_{n\nu }$ and $d_{n\nu }$ are the ramification index and the
exponent of the different for the extension $K_\nu /K_n$. On the
other hand [W], $d_\nu =e_{n\nu }d_n+d_{n\nu }$, so that $l=d_n$,
and we obtain that $x\in S_n$.

Hence, $S_\nu \subset S_n$ for $\nu >n$. It follows from (8) that
$\mu (S_n)=1$ for each $n$. Since $S=\bigcap \limits _{n=1}^\infty S_n$,
we get the desired equality $\mu (S)=1$.

(b) It follows from the definition of the measure $\mu $ by
the relation (7) that we need only to prove, for any $y\in S$ and
any $a\in K$ with $\|a\|\le 1$, the inequality
\begin{equation}
|\langle a,y\rangle |_1\le 1.
\end{equation}

Suppose that $a\in K_n,\ \|a\|\le 1$. Then
$\langle a,y\rangle =T(ay_n)=\mbox{Tr}_{K_n/k}(m_n^{-1}ay_n)$, $\ \ y_n
=T_n(y)$.
Since $\|y_n\|\le \|m_n\|q_n^{d_n/m_n}$ we find that
\[
|m_n^{-1}ay_n|_n\le q_n^{d_n}
\]
so that the inequality (11) follows from the definition of $d_n$.

(c) This property is an immediate consequence of (a) and the ultrametric
property of absolute values.

(d) We have shown (while proving Proposition 4) that
\[
\mu \left( \left\{ x\in \overline{K}:\ |\langle a,x\rangle |_1>
q_1^{-[-N/e_n]}\right\} \right) =0.
\]
This implies the inequality
\[
\mbox{ess sup }_{x\in
\overline{K}}\ |\langle a,x\rangle |_1\le q_1^{-[-N/e_n]}
\]

To prove the inverse inequality recall that $\langle a,x\rangle
=m_n^{-1}\mbox{Tr}_{K_n/k}(aT_n(x))$. Suppose that
$|m_n|_n=q_n^\nu,\ \nu \in \mbox{{\bf Z}}$. By properties of the trace [W]
there exists such an element $y_n\in K_n$ that $|y_n|_n=q_n^{\nu +d_n}$,
\[
|\mbox{Tr}_{K_n/k}(ay_n)|_1=q_1^l
\]
where $l\in \mbox{{\bf Z}},\ e_n(l-1)<N+\nu \le e_nl$, whence
\[
|T(ay_n)|_1=q_1^{l-\nu /e_n}\ge q_1^{-[-N/e_n]}.
\]
If a natural number M is sufficiently large then the inequality \linebreak
$\|z-y_n\|\le q_1^{-M} \ \ (z\in K_n)$ implies the equality
$|T(az)|_1=q_1^{l-\nu /e_n}$.

Consider a subset in $\overline{K}$ of the form
\[
U=\left\{x\in \overline{K}:\ \ \|T_n(x)-y_n\|\le q_1^{-M}\right\},
\]
with the number $M$ mentioned above. If $x\in U$ then
\[
|\langle a,x\rangle |_1=|T(aT_n(x))|_1=q_1^{l-\nu /e_n}\ge q_1^{-[-N/e_n]}.
\]
On the other hand, using the formula (8) and taking into account
that $\|y_n\|=\|m_n\|q_n^{d_n/m_n}$ we find that $\mu (U)\ne 0$
if $M$ is taken sufficiently large. This means that
\[
\mbox{ess sup }_{x\in \overline{K}}\ |
\langle a,x\rangle |_1\ge q_1^{-[-N/e_n]}
\]

(e) Consider subsets in $\overline{K}$ having the form
$A_{n,B}=\{x\in \overline{K}:\ \ T_n(x)\in B\}$ where $B\subset K_n$
is a Borel set, $n\ge 1$. If $\nu >n$ then
$A_{n,B}=\{x\in \overline{K}:\ \ T_\nu (x)\in C\}$
where $C=\{y\in K_\nu :\ \ T_n(y)\in B\}$. It follows from this
observation that the sets $A_{n,B}$ form a Boolean algebra of sets
which of course generates ${\cal B}(\overline{K})$.

By general principles of measure theory this implies the existence,
for any $V\in {\cal B}(\overline{K})$ and any $\varepsilon >0$, of
such a set $A_{n,B}$ that $\mu (V\Delta A_{n,B})<\varepsilon $
where $\Delta $ means symmetric difference. Here one can
take $B\subset \{y\in K_n:\ \ \|y\|\le q_n^{d_n/m_n}\|m_n\|\}$.

Hence an indicator of any set from ${\cal B}(\overline{K})$ can
be approximated in $L_s(\overline{K})$ by indicators of the sets
$A_{n,B}$ which belong to ${\cal D}(\overline{K}).\ \ \Box $

Let ${\cal E}(\overline{K})$ be the set of functions
$f(x)=\varphi (T_n(x))$ where $n\ge 1,\ \varphi $ is a locally
constant function on $K_n$. It follows from Theorem 1 that
${\cal E}(\overline{K})\subset L_1(\overline{K})$, and the
formula (8) remains valid for $f\in {\cal E}(\overline{K})$.

4. In the case of Galois extensions (see Sect. 2.2) the measure
$\mu $ is invariant with respect to the Galois group $G$.
This follows immediately from (1) and (7).

\section{Fourier Transform and Fractional Differentiation}

1. We define the Fourier transform of a function
$f\in L_1(\overline{K})$ as
\[
\widehat{f}(\xi )=\int _{\overline{K}}\chi (\langle \xi ,x\rangle )
f(x)\,d\mu (x),\ \ \xi \in K.
\]
Note that to avoid confusion we use here the notation $\widehat{f}$
in contrast to the notation $\widetilde{\varphi }$
for the Fourier transform over a local field. Sometimes we shall
also write ${\cal F}f=\widehat{f}$.

Let $f\in {\cal E}(\overline{K}),\ f(x)=\varphi (T_n(x))$. An
explicit form of $\widehat{f}(\xi )$ can be easily obtained from (8):
if $\xi \in K_\nu ,\ \nu \ge n$, then
\begin{equation}
\widehat{f}(\xi )=q_\nu ^{-d_\nu }\|m_\nu \|^{-m_\nu }
\int \limits_{\zeta \in K_\nu :\ \|\zeta \|\le q_\nu ^{d_\nu /m_\nu }
\|m_\nu \|}\chi (T(\xi \zeta ))\varphi (T_n(\zeta ))\,d\zeta .
\end{equation}
Writing (12) for $\nu =n$ and using the inversion formula (5) we find that
\begin{equation}
\varphi (z)=\int _{K_n}\chi (-T(\xi z))\widehat{f}(\xi )\,d\xi ,
\ \ z\in K_n,\ \|z\|\le q_n^{d_n/m_n}\|m_n\|.
\end{equation}

Let us describe the image of ${\cal E}(\overline{K})$ under the
transform ${\cal F}$. For $\xi \in K$ we shall write
\[
\mbox{dist}(\xi ,K_n)=\inf \limits _{\eta \in K_n}\|\xi -\eta \|.
\]

\begin{teo}
A function $F(\xi ),\ \xi \in K$, is an image of some function
$f\in {\cal E}(\overline{K}),\ f(x)=\varphi (T_n(x))$, under
the transform ${\cal F}$, if and only if: (i) the restriction of
$F$ to each field $K_\nu ,\ \nu \ge n$, belongs to ${\cal D}(K_\nu )$;
(ii) $F(\xi +\xi ')=F(\xi )$ for any $\xi ,\xi '\in K$ with
$\|\xi '\|\le 1$; (iii) $F(\xi )=0$ if {\em dist}$(\xi ,K_n)>1$.
\end{teo}

{\it Proof.}  The necessity of (i) and (ii) follows immediately
from (12). According to (13),
\begin{equation}
\varphi (T_n(\zeta ))=
\int _{K_n}\chi (-T(\zeta \eta ))\widehat{f}(\eta )\,d\eta ,
\end{equation}
if $\zeta \in K_\nu ,\ \nu >n,\ \|T_n(\zeta )\|\le q_n^{d_n /m_n}\|m_n\|$.
As we have seen, the last inequality is valid if
$\|\zeta \|\le q_\nu ^{d_\nu /m_\nu }\|m_\nu \|$. Thus we may
substitute (14) into (12). We obtain for any $\xi \in K_\nu $ that
\[
F(\xi )=\widehat{f}(\xi )
=q_\nu ^{-d_\nu }\|m_\nu \|^{-m_\nu }\int _{K_n}F(\eta )\,d\eta
\int \limits_{\zeta \in K_\nu :\ \|\zeta \|\le q_\nu ^{d_\nu /m_\nu }
\|m_\nu \|}\chi (T(\zeta (\xi -\eta )))\,d\zeta =
\]
\[
=\int \limits_{\eta \in K_n:\ \|\eta -\xi \|\le 1}F(\eta )\,d\eta
\]
by virtue of (2). Using (ii) we find that
\[
F(\xi )=F(\xi )\int _{\eta \in K_n:\ \|\eta -\xi \|\le 1}d\eta .
\]

If $\mbox{dist}(\xi ,K_n)>1$ then the domain of integration is
empty so that $F(\xi )=0$.

Now we turn to the sufficiency of (i)-(iii). Put
\begin{equation}
\varphi (z)=\int _{K_n}\chi (T(-\xi z))F(\xi )\,d\xi,\ \ z\in K_n,
\end{equation}
$f(x)=\varphi (T_n(x)),\ \ x\in \overline{K}$.
Let us compute $\widehat{f}$. Suppose that
$\xi \in K_\nu ,\ \nu \ge n$. Then we may substitute (15) into
(12). Repeating the previous calculation we obtain that
\[
\widehat{f}(\xi )=F(\xi )\int _{\eta \in K_n:\ \|\eta -\xi \|\le 1}d\eta .
\]

If $\mbox{dist}(\xi ,K_n)>1$ then $\widehat{f}(\xi )=0=F(\xi )$.
Let $\mbox{dist}(\xi ,K_n)\le 1$. There exists such an element
$\sigma \in K_n$ that $\mbox{dist}(\xi ,K_n)=\|\xi -\sigma \| $
(see [NBB]). Thus $\|\xi -\sigma \|\le 1$, and by the
ultra-metric inequality
\[
\{\eta \in K_n:\ \|\eta -\xi \|\le 1\}
=\{\eta \in K_n:\ \|\eta -\sigma \|\le 1\}.
\]
whence
\[
\int _{\eta \in K_n:\ \|\eta -\xi \|\le 1}d\eta
=\int _{\|\eta -\sigma \|\le 1}d\eta =1
\]
due to the invariance of the Haar measure on $K_n$. Hence
$\widehat{f}(\xi )=F(\xi )$ for this case too.$\ \ \ \Box $

2. Using (8), (12) and the isometric property of the Fourier
transform over a local field we come to the Plancherel formula
for $f(x)=\varphi (T_n(x))\ , f\in {\cal E}(\overline{K})$:
if $\nu \ge n$ then
\[
\int _{\overline{K}}|f(x)|^2\,d\mu (x)
=\int _{K_\nu }|\widehat{f}(\xi )|^2\,d\xi=
\]
\[
=q_n^{-d_n}\|m_n\|^{-m_n}\int \limits_{z\in K_n:\ \|z\|
\le q_n^{d_n/m_n}\|m_n\|}|\varphi (z)|^2\,dz.
\]

If $g(x)=\psi (T_m(x))$ is another function
from ${\cal E}(\overline{K})$ and $\nu \ge \max (m,n)$ then
\[
\int _{\overline{K}}f(x)\overline{g(x)}\,d\mu (x)
=\int _{K_\nu }\widehat{f}(\xi )\overline{\widehat{g}(\xi )}\,d\xi .
\]

Hence, the transform ${\cal F}$ is the isometry of
$L_2(\overline{K})$ onto the completion of the image
${\cal F}{\cal E}(\overline{K})$ (described in Theorem 2) with respect
to the norm
\[
\lim \limits _{\nu \to \infty }\int _{K_\nu }|F(\xi )|^2\,d\xi .
\]

3. Next we shall study an operator $D^\alpha $ defined (initially)
on ${\cal E}(\overline{K})$ as
\[
D^\alpha ={\cal F}^{-1}\Delta ^\alpha {\cal F},\ \ \alpha >0,
\]
where $\Delta ^\alpha $ is an operator of multiplication in
${\cal F}{\cal E}(\overline{K})$ by the function
\[
\Delta ^\alpha (\xi )=\left\{ \begin{array}{rl}
\|\xi \|^\alpha , & \mbox{if }\|\xi \|>1\\
0, & \mbox{if }\|\xi \|\le 1
\end{array} \right. ,\ \ \xi \in K.
\]
Theorem 2 shows that $D^\alpha $ is well-defined. It is easy to
see that  $D^\alpha $ is essentially self-adjoint as an operator
on $L_2(\overline{K})$.

Let us begin with a simple property of $D^\alpha $ justifying
its interpretation as an analogue of fractional differentiation.

Let $f_a(x)=\chi (\langle a,x\rangle ),\ x\in \overline{K}$, where
$a\in K$. Note that if $\|a\|\le 1$ then $f_a(x)=1$ for all $x\in S$
or, in other words, $f_a(x)=1$ $\mu $-almost everywhere on $\overline{K}$.

\begin{prop}
For $\mu $-almost all $x\in \overline{K}$
\[
(D^\alpha f_a)(x)=\left\{ \begin{array}{rl}
0, & \mbox{if }\|a\|\le 1\\
\|a\|^\alpha f_a(x), & \mbox{if }\|a\|> 1
\end{array} \right.
\]
\end{prop}

{\it Proof.} Suppose that $a\in K_n$. Then $f_a(x)=\varphi (T_n(x))$
where $\varphi (z)=\chi (T(az)),\ z\in K_n$. Using (12) and (2)
we find that
\[
\widehat{f_a}(\xi )=\Omega _n(\xi +a),\ \ \xi \in K_n,
\]
where
\[
\Omega _n(z)=\left\{ \begin{array}{rl}
1, & \mbox{if }\|z\|\le 1\\
0, & \mbox{if }\|z\|> 1
\end{array} \right.
\]

Thus $(D^\alpha f_a)(x)=\psi (T_n(x))$ where
\[
\psi (z)=\int _{K_n}\chi (T(-z\xi ))\Delta ^\alpha (\xi )
\Omega _n(\xi +a)\,d\xi =
\]
\[
=\int \limits_{\|\xi \|>1,\ \|\xi +a\|\le 1}\|\xi \|^\alpha \chi (T(-z\xi ))\,d\xi
\]
for $\|z\|\le q_n^{d_n/m_n}\|m_n\|$. If $\|a\|\le 1$ then
$\psi (z)=0$. If $\|a\|>1$ then $\|\xi \|=\|a\|$ on the whole
integration domain so that
\[
\psi (z)=\|a\|^\alpha \int _{\|\xi +a\|\le 1}\chi (T(-z\xi ))\,d\xi
=\|a\|^\alpha \chi (T(az))
\]
if  $\|z\|\le q_n^{d_n/m_n}\|m_n\|.\ \ \ \Box $

It follows from the definition of $D^\alpha $ that
\[
\int _{\overline{K}}(D^\alpha f)(x)\,d\mu (x)=0
\]
for any $f\in {\cal E}(\overline{K})$. In other words, $\mu $ is
a {\it harmonic measure} for $D^\alpha $.

Though we are interested primarily in the operator $D^\alpha $
with $\alpha >0$ (in view of some applications presented in
Sect. 5 below), we shall need also the case $\alpha \le 0$.
An elementary calculation shows that
\begin{equation}
(D^0f)(x)=f(x)-\int _{\overline{K}}f(y)\,d\mu (y),\ \ x\in S,
\end{equation}
for any $f\in {\cal E}(\overline{K})$.

It follows from Theorem 2 that $(D^\alpha f)(x),\ x\in S$, is
well-defined (if $f\in {\cal E}(\overline{K})$) for all
$\alpha \in \mbox{{\bf C}}$ being in fact an entire function with
respect to $\alpha $. It is clear that
$D^{\alpha _1}D^{\alpha _2}f
=D^{\alpha _1+\alpha _2}f,\ \ \alpha _1,\alpha _2\in \mbox{{\bf C}},
\ \ f\in {\cal E}(\overline{K})$.

\begin{prop}
If $f\in {\cal E}(\overline{K}),\int \limits _{\overline{K}}f(x)\,d\mu (x)
=0$ then the general solution $u\in
{\cal E}(\overline{K})$ of the equation
\begin{equation}
D^\alpha u=f,\ \ \ \alpha >0,
\end{equation}
is of the form $u=D^{-\alpha }f+C$ where $C$ is an arbitrary constant.
\end{prop}

{\it Proof.}  The function $u$ is a solution of (17) in view
of the relation (16), since $D^\alpha C=0$ by Proposition 5.
On the other hand, let $v\in {\cal E}(\overline{K}),\ D^\alpha v=0$
$\mu $-almost everywhere on $\overline{K}$. By the definition
of $D^\alpha $, we find that $\widehat{v}(\xi )=0$ if $\|\xi \|>1$.
Now if $v(x)=\varphi (T_n(x))$ then
\[
\varphi (z)=\int _{\|\xi \|\le 1}\chi (T(-z\xi ))\widehat{v}(\xi )\,d\xi
=\widehat{v}(0)
\]
for all $z\in K_n$ with $\|z\|\le q_n^{d_n/m_n}\|m_n\|$. This means
that $v(x)=\mbox{const}$ $\mu $-almost everywhere on
$\overline{K}.\ \ \ \Box $

4. The fractional differentiation operator on a local field admits
a hyper-singular integral representation not containing the
Fourier transform [VVZ]. Now we give a representation of this kind
for our present situation.

\begin{prop}
Let $\alpha >0,\ f\in {\cal E}(\overline{K}),
\ f(x)=\varphi (T_n(x))$. Then $(D^\alpha f)(x)=\psi (T_n(x))$ where
\[
\psi (z)=q_n^{d_n\alpha /m_n}\frac{1-q_n^{\alpha /m_n}}{1-q_n^{-1-
\alpha /m_n}}\|m_n\|^{-m_n}\int \limits_{x\in K_n:\  \|x\|
\le q_n^{d_n/m_n}\|m_n\|}\biggl[ \|x\|^{-m_n-\alpha }\|m_n\|^{m_n+
\alpha }\biggr.
\]
\begin{equation}
\left. +\frac{1-q_n^{-1}}{q_n^{\alpha /m_n}-1}q_n^{-d_n(1+\alpha /m_n)}
\right][\varphi (z-x)-\varphi (z)]\,dx,
\end{equation}
$z\in K_n,\ \ \|z\|\le q_n^{d_n/m_n}\|m_n\|$.
\end{prop}

{\it Proof.}  Let us compute the Fourier transform over $K_n$ for
the function $\Delta ^\lambda (\xi )$ with $\Re \lambda <-1$. We have
\[
\widetilde{\Delta ^\lambda }(x)
=q_n^{-d_n/2}\sum \limits _{N=1}^\infty q_n^{N
\lambda /m_n}\int _{|\xi |_n=q_n^N}\left( \chi
\circ \mbox{Tr}_{K_n/k}\right) (x\xi )\,d\xi .
\]
Using (3) and summing the geometric progression we obtain after
elementary transformation that
\[
\widetilde{\Delta ^\lambda }(x)=q_n^{d_n(1/2+\lambda /m_n)}
\frac{1-q_n^{\lambda /m_n}}{1-q_n^{-1-\lambda /m_n}}\left(
|x|_n^{-1-\lambda /m_n}+\right.
\]
\begin{equation}
\left. +q_n^{-d_n(1+\lambda /m_n)}\frac{1-q_n^{-1}}{q_n^{\lambda /m_n}-1}
\right) \omega _n(x),\ \ x\in K_n,
\end{equation}
where
\[
\omega _n(x)=\left\{ \begin{array}{rl}
1, & \mbox{if }|x|_n\le q_n^{d_n}\\
0, & \mbox{if }|x|_n>q_n^{d_n}
\end{array} \right. .
\]

Under our normalizations the relation between Fourier transform
and convolution reads
$\widetilde{u\ast v}=q_n^{d_n/2}\widetilde{u}\widetilde{v}$ (with
the usual assumptions with regard to the functions $u,v$). Thus
for $\Re \lambda <-1$ we get $(D^\lambda f)(x)=\psi _\lambda (T_n(x))$
where
\[
\psi _\lambda (m_ny)=\left( \widetilde{\Delta ^\lambda }\ast
\varphi _n\right) (y),\ \ y\in K_n,\ \ \|y\|\le q_n^{d_n/m_n},
\]
$\widetilde{\Delta ^\lambda }$ is given by (19),
\[
\varphi _n(\eta )=\left\{ \begin{array}{rl}
q_n^{-d_n/2}\varphi (m_n\eta ), & \mbox{if }|\eta |_n\le q_n^{d_n}\\
0, & \mbox{if }|\eta |_n>q_n^{d_n}
\end{array} \right. .
\]

Since $\Delta ^\lambda (0)=0$ we find that
\[
\int _{|x|_n\le q_n^{d_n}}\widetilde{\Delta ^\lambda }(x)\,dx
=0,\ \ \ \Re \lambda <-1.
\]
Thus the above convolution may be written as follows:
\[
\psi _\lambda (m_ny)=q_n^{d_n\lambda /m_n}
\frac{1-q_n^{\lambda /m_n}}{1-q_n^{-1-\lambda /m_n}}
\int _{|x|_n\le q_n^{d_n}}\biggl[ |x|_n^{-1-\lambda /m_n}+\biggr.
\]
\begin{equation}
\left. +\frac{1-q_n^{-1}}{q_n^{\lambda /m_n}-1}q_n^{-d_n(1+
\lambda /m_n)}\right][\varphi (m_ny-m_nx)-\varphi (m_ny)]\,dx
\end{equation}

Since the left-hand side of (20) is an entire function with
respect to $\lambda $, the analytic continuation of the
right-hand side leads to the required equality (18). \ \ \ $\Box $

\section{Heat Equation}

1. We begin with an auxiliary equation over the local field $K_n$.
Note that this equation is different from a local field heat equation
studied in [Ha1, K1, VVZ].

Let $\partial _n^\alpha $ be a pseudo-differential operator over $K_n$
with the symbol $\Delta ^\alpha (\xi )$, $\ \alpha >0$. It means that for any
$\varphi \in {\cal D}(K_n)$ \ $(\partial _n^\alpha \varphi )(x),\ x\in
K_n$, is the inverse Fourier transform of the function $\Delta ^\alpha (\xi )
\widetilde{\varphi }(\xi ),\ \xi \in K_n$. As we have proved (in the course
of proving Proposition 7)
\[
(\partial _n^\alpha \varphi )(x)=q_n^{d_n\alpha /m_n}\frac{1-q_n^{\alpha
/m_n}}{1-q_n^{-1-\alpha /m_n}}\int _{|x|_n\le q_n^{d_n}}\biggl[ |x|_n^{-1
-\alpha/m_n}+\biggr.
\]
\begin{equation}
\left. +q_n^{-d_n(1+\alpha /m_n)}\frac{1-q_n^{-1}}{q_n^{\alpha /m_n}-1}
\right] [\varphi (x-y)-\varphi (x)]\,dy,
\end{equation}
$\varphi \in {\cal D}(K_n)$.

Consider the equation
\begin{equation}
\frac{\partial u(x,t)}{\partial t}+\left( \partial _n^\alpha u\right)
(x,t)=0,\ \ \ x\in K_n,\ t>0,
\end{equation}
with an initial condition $u(x,0)=\varphi (x),\ \ \varphi \in {\cal D}(K_n)$.
As usual, a fundamental solution of this Cauchy problem is
\begin{equation}
\Gamma _\alpha ^{(n)}(x,t)=q_n^{-d_n}\int _{K_n}(\chi \circ \mbox{Tr}_{K_n
/k})(-x\xi )\rho _\alpha (\|\xi \|,t)\,d\xi
\end{equation}
where
\begin{equation}
\rho _\alpha (s,t)=\left\{ \begin{array}{rl}
e^{-ts^\alpha }, & \mbox{if } s>1 \\
1, & \mbox{if } s\le 1,
\end{array} \right.
\end{equation}
$s\ge 0,\ t>0$ (note the connection between convolution and the Fourier
transform; see Sect. 4.4).

It is clear that $\Gamma _\alpha ^{(n)}$ is continuous on $K_n\times
(0,\infty )$. If $\|x\|>q_n^{d_n/m_n}$ then we can find such an element
$\xi _0\in K_n,\ \|\xi _0\|=1$, that $(\chi \circ \mbox{Tr}_{K_n/k})
(-x\xi _0)\ne 1$. After the change $\xi =\xi _0+\eta $ in (23) we find that
$\Gamma _\alpha ^{(n)}(x,t)=0$. This observation allows to differentiate
under the sign of integral showing that the function
\begin{equation}
u(x,t)=\int _{K_n}\Gamma _\alpha ^{(n)}(x-y,t)\varphi (y)\,dy
\end{equation}
satisfies (22).

Since $\int \limits_{K_n}\Gamma _\alpha ^{(n)}(x,t)\,dx=1$ we may write
\begin{equation}
u(x,t)=\varphi (x)+\int \limits_{\|y-x\|\le q_n^{d_n/m_n}}\Gamma _\alpha ^{(n)}
(x-y,t)[\varphi (y)-\varphi (x)]\,dy.
\end{equation}

On the other hand, let $|x|_n=q_n^N,\ N\le d_n$. It follows from (23),(24)
and (3) that
\[
q_n^{d_n}\Gamma _\alpha ^{(n)}(x,t)=1+(1-q_n^{-1})\sum \limits _{j=1}^{d_n
-N}q_n^j\rho _\alpha (q_n^j,t)-\rho _\alpha (q_n^{d_n-N+1},t)q_n^{d_n-N}.
\]
This representation shows that $\Gamma _\alpha ^{(n)}(-x,t)
=\Gamma _\alpha ^{(n)}(x,t)$, $\Gamma _\alpha ^{(n)}(x,t)\to 0$ when
$x\ne 0$ and $t\to 0$,
\[
|\Gamma _\alpha ^{(n)}(x,t)|\le 2q_n^{d_n/2}|x|_n^{-1},\ \ \ x\ne 0.
\]
Returning to (26), using the local constancy of $\varphi $ and the
dominated convergence theorem we find that $u(x,t)\to \varphi (x)$ as
$t\to 0$.

Finally, let us prove that $\Gamma _\alpha ^{(n)}(x,t)\ge 0$ for all
$x,t$. It is sufficient to show that the function (25) is
non-negative for any non-negative $\varphi \in {\cal D}(K_n)$.

Let us fix $T>0$. Suppose that, on the contrary, $u(x,t)<0$ for some
values of $x\in K_n,\ t\in [0,T]$. Since $u(x,t)$ vanishes as a
function of $x$ outside a compact set not depending on $t\in [0,T]$,
and $\lim \limits _{t\to 0}u(x,t)\ge 0$, the function $u(x,t)$ attains
its negative global minimum at a certain point $(x_0,t_0)$ where
$t_0>0$.

It follows from (21) that $\left( \partial _n^\alpha u\right) (x_0,
t_0)\le 0$. Of course, also
${\displaystyle \frac{\partial u}{\partial t}}(x_0,t_0)\le 0$
so that by (22) we get
\[
\frac{\partial u}{\partial t}(x_0,t_0)=0,\ \
\left( \partial _n^\alpha u\right) (x_0,t_0)=0.
\]
The second equality means (see (21)) that $u(x_0,t_0)$ does not depend
on $x$. Since $u(x,t_0)$ has a compact support, it follows that
$u(x_0,t_0)=0$, and we have come to a contradiction.

2. The relation (23) and the above properties of $\Gamma _\alpha ^{(n)}$
lead to the following result.

Let us consider $\rho _\alpha (\|\xi \|,t)$ as a function on $K\times
(0,\infty )$.

\begin{prop}
For each $t>0$ the function $\rho _\alpha (\|\xi \|,t)$ is a continuous
positive definite function on $K$.
\end{prop}

Now we are again in a position to use Madrecki's theorem [Ma], thus
obtaining, for each $t>0$, a Radon probability measure $\pi (t,dx)$ on
${\cal B}(\overline{K})$ such that
\begin{equation}
\rho _\alpha (\|\xi \|,t)=\int _{\overline{K}}\chi (\langle \xi ,x
\rangle )\pi (t,dx),\ \ \ \xi \in K.
\end{equation}

Just as in Sect. 3, we show that if $f(x)=\varphi (T_n(x))$, $\varphi
\in {\cal D}(K_n)$, then
\begin{equation}
\int _{\overline{K}}f(x)\pi (t,dx)
=\|m_n\|^{-m_n}\int \limits_{z\in K_n:\ \|z\|
\le q_n^{d_n/m_n}\|m_n\|}\Gamma _\alpha ^{(n)}(m_n^{-1}z,t)\varphi (z)\,dz.
\end{equation}
It follows from (28) that the measures $\pi (t,dx)$ are concentrated
on the set $S$, and (28) remains valid for any
$f\in {\cal E}(\overline{K})$.

Since $\rho _\alpha (\|\xi \|,t_1+t_2)=\rho _\alpha (\xi ,t_1)\rho _\alpha
(\xi ,t_2)$, we find that $\pi (t_1+t_2,\cdot )=\pi (t_1,\cdot )\ast \pi
(t_2,\cdot )$. Hence the family of measures $\pi (t,dx)$ determines a
spatially homogeneous Markov process $X(t)$ on $\overline{K}$ with the
transition probabilities $P(t,x,\Lambda )=\pi (t,\Lambda -x)$, $\Lambda
\in {\cal B}(\overline{K})$.

Let $U_t$ be a corresponding semigroup of operators on the Banach space of
bounded measurable functions on $\overline{K}$:
\[
(U_tf)(x)=\int _{\overline{K}}f(x+y)\pi (t,dy).
\]

\begin{teo}
(a) $X(t)$ is a stochastically continuous and right-continuous \newline
Markov process on $\overline{K}$ without discontinuities of the second
kind.\newline
(b) The operator $-D^\alpha $ is the generator of the process $X(t)$,
in the sense that for any $f\in {\cal E}(\overline{K})$ the function
$u(x,t)=(U_tf)(x)$ satisfies the equation
\[
\frac{\partial u}{\partial t}+D^\alpha u=0.
\]
(c) The measure $\mu $ is an invariant measure for the process $X(t)$.
\end{teo}

{\it Proof}. Let $f(x)=\varphi (T_n(x)),\ \ \varphi \in {\cal E}(K_n)$.
In accordance with (28), $u(x,t)=F_t(T_n(x))$ where
\[
F_t(\zeta )=\|m_n\|^{-m_n}\int \limits_{\|z\|\le q_n^{d_n/m_n}\|m_n\|}\Gamma
_\alpha ^{(n)}(m_n^{-1}z,t)\varphi (\zeta +z)\,dz,
\]
$\zeta \in K_n$. By the results from the preceding section,
${\displaystyle \frac{\partial
u}{\partial t}}+D^\alpha u=0$ and also $u(x,t)\to f(x)$ as $t\to 0$
which implies stochastic continuity.

Let us show that $X(t)$ is a cadlag process. By the Kinney-Dynkin
criterion [Dy] and the definition of the topology in $\overline{K}$ it is
sufficient to prove that for any compact set $B\subset \overline{K}$,
any $N=1,2,\ldots $, and any $\varepsilon >0$
\begin{equation}
\lim \limits _{t\to 0}\sup \limits _{x\in B}P(t,x,C_{\varepsilon ,N}(x))
=0
\end{equation}
where
\[
C_{\varepsilon ,N}(x)=\bigcup \limits _{n=1}^N\{y\in \overline{K}:\ \
\|T_n(x-y)\|\ge \varepsilon \}.
\]
We find by virtue of (28) that
\[
P(t,x,C_{\varepsilon ,N}(x))\le \sum \limits _{n=1}^N\int \limits
_{\|T_n(y)\|\ge \varepsilon }\pi (t,dy)=
\]
\[
=\sum \limits _{n=1}^N\|m_n\|^{-m_n}\int \limits _{z\in K_n:\
\varepsilon \le \|z\|\le q_n^{d_n/m_n}\|m_n\|}\Gamma _\alpha ^{(n)}
(m_n^{-1}z,t)\,dz,
\]
and the equality (29) follows from the properties of
$\Gamma _\alpha ^{(n)}$ established above.

Finally, let
\[
\mu _1(A)=\int _{\overline{K}}P(t,x,A)\,d\mu (x),\ \ \ A\in {\cal B}(
\overline{K}),
\]
and we need to prove that $\mu _1=\mu $. Indeed, for any $\xi \in K$
\[
\int _{\overline{K}}\chi (\langle \xi ,x\rangle )\,d\mu _1(x)=\int
_{\overline{K}}d\mu (y)\int _{\overline{K}}\chi (\langle \xi ,x\rangle
)P(t,y,dx)=
\]
\[
=\int _{\overline{K}}\chi (\langle \xi ,y\rangle )\,d\mu (y)
\int _{\overline{K}}\chi (\langle \xi ,x\rangle )\pi (t,dx)=
\rho _\alpha (\|\xi \|,t)\Omega (\xi )=\Omega (\xi )
\]
whence $\mu _1=\mu .\ \ \ \ \Box $

{\it Remark. } If $x\in S,\ \Lambda \cap S=\emptyset$, then
$(\Lambda -x)\cap S=\emptyset $. Indeed, if $y\in \Lambda $
then
$$
\Vert T_n(y)\Vert >q_n^{d_n/m_n}\Vert m_n\Vert
$$
whence $\Vert T_n(y-x)\Vert =\Vert T_n(y)\Vert $ due to the
ultra-metric property, so that $y-x\not\in S$. Hence, a
trajectory of the process $X$ starting at a point from $S$
remains in $S$ almost surely.

This means that the part $X_S$ of the process $X$ in $S$ is a
non-exploding process. Since $S$ is compact in a strong
topology of $\overline{K}$ (which is a consequence of general
properties of projective limits [RR]), it follows from Theorem
I.9.4 of [BG] that $X_S$ is a Hunt process.

3. If we deal with Galois extensions then as before we see from the
definition of $P$ that this transition probability is invariant with
respect to the Galois group $G$.

\section* {Acknowledgment}

This research was supported in part by grants from the International
Science Foundation and the Ukrainian Fund for Fundamental Research.

\end{document}